\magnification=1200 \baselineskip=18pt
 \vsize=8.7truein
\hsize=5.5 truein \overfullrule0pt \hoffset=0.35in \voffset=0.3in
\overfullrule0pt
\parskip=4pt
\parindent=20pt
\nopagenumbers \pageno=1 \font\ninerm=cmr9
\headline={\ifnum\pageno>1\hfill\ninerm\folio\else\fi}
\def\foot#1#2{{\baselineskip=9pt\sevenrm\footnote{$^{#1}$}{#2}}}
\def\nh{\noindent\hangindent=1.5truecm\hangafter=1}
\def\bs{\bigskip}
\def\cl{\centerline}

\def\ni{\noindent}
\def\ve{\vfill\eject}
\def\vend{\ve\end}

\def\a{\alpha}
\def\p{p}
\def\ab{\allowbreak}
\def\as{^*}

\def\be{\beta}

\def\bigmi{\;\big|\;}

\def\bY{{\bar Y}}

\def\cE{{\cal E}}

\def\cP{{\cal P}}

\def\cY{{\cal Y}}

\def\dag{^{\dagger}}
\def\de{\delta}

\def\ep{\epsilon}
\def\etal{{\it et al.}~}
\def\ga{\gamma}

\def\half{^{1/2}}

\def\hf{{\hat f}}
\def\hF{{\widehat F}}

\def\hka{{\hat\ka}}

\def\hP{{\hat P}}

\def\hatt{{\hat t}}
\def\hth{{\hat\th}}

\def\ia{_{i\a}}
\def\ij{_{ij}}

\def\ka{\kappa}
\def\la{\lambda}

\def\mhf{^{-1/2}}
\def\mi{\;|\;}
\def\mo{^{-1}}

\def\oon{{1\over n}}

\def\oi{_{1i}}

\def\otd{{\textstyle{1\over3}}}

\def\ra{\rightarrow}
\def\rai{\ra\infty}

\def\sumi{\sum_i\,}

\def\sumionu{\sum_{i=1}^\nu\,}

\def\sumj{\sum_j\,}
\def\sumjon{\sum_{j=1}^n\,}
\def\sumk{\sum_k\,}

\def\th{\theta}
\def\tthd{{\textstyle{2\over3}}}
\def\thf{{\textstyle{1\over2}}}

\def\var{{\rm var}}

\def\zi{_{0i}}
\def\nu{N}

\cl{\bf TO HOW MANY SIMULTANEOUS HYPOTHESIS TESTS}
\cl{\bf CAN NORMAL, STUDENT'S $t$}
\cl{\bf OR BOOTSTRAP CALIBRATION BE APPLIED?}

\bs

\centerline{Jianqing Fan \quad Peter Hall \quad Qiwei Yao}

\foot{}{Jianqing Fan is Frederick L. Moore '18 Professor of Finance
(E-mail: jqfan@princeton.edu), Department of Operations Research and
Financial Engineering, Princeton University, Princeton, NJ 08544,
USA. Peter Hall (E-mail: halpstat@maths.anu.edu.au) is professor,
Centre for Mathematics and its Applications, Australian National
University, Canberra, ACT 0200, Australia. Qiwei Yao is professor
(E-mail: q.yao@lse.ac.uk), Department of Statistics, London School
of Economics, Houghton Street, London WC2A 2AE, United Kingdom.
Fan's work was sponsored in part by NSF grant DMS-0354223 and NIH
grant R01-GM07261.}

\bs\bs

\baselineskip=12.5pt

\noindent{\bf ABSTRACT.}  In the analysis of microarray data, and in
some other contemporary statistical problems, it is not uncommon to
apply hypothesis tests in a highly simultaneous way.  The number,
$\nu$ say, of tests used can be much larger than the sample sizes,
$n$, to which the tests are applied, yet we wish to calibrate the
tests so that the overall level of the simultaneous test is
accurate.  Often the sampling distribution is quite different for
each test, so there may not be an opportunity for combining data
across samples.  In this setting, how large can $\nu$ be, as a
function of $n$, before level accuracy becomes poor?  In the present
paper we answer this question in cases where the statistic under
test is of Student's $t$ type.  We show that if either Normal or
Student's $t$ distribution is used for calibration then the level of
the simultaneous test is accurate provided $\log\nu$ increases at a
strictly slower rate than $n^{1/3}$ as $n$ diverges.  On the other hand, if bootstrap methods are used for calibration then we may choose $\log\nu$ almost as large as $n\half$ and still achieve asymptotic level accuracy.  The implications of these results are explored both theoretically and numerically.

\bs\bs







\baselineskip=20pt

\eject

\cl{\bf 1.  INTRODUCTION}

Modern technology allows us to collect a large amount
of data in one scan of images.  This is exemplified in genomic
studies using microarrays, tiling arrays and proteomic techniques.
In the analysis of microarray data, and in some other contemporary
statistical problems, we often wish  to make statistical inference
simultaneously for all important parameters. The number of
parameters, $\nu$, is frequently much larger than sample size,~$n$.  Indeed,
sample size is typically small; e.g.~$n=8$, 20 or 50
are considered to be typical, moderately large or large, respectively, for
microarray data. The question arises naturally as to how large $\nu$ can
be before the accuracy of simultaneous statistical inference becomes poor.

Important results in this direction have been obtained by van der Laan and Bryan (2001), who showed that the population mean and variance parameters can be consistently estimated when $\log\nu=o(n)$ if observed data are bounded.  Bickel and Levina (2004) gave similar results in a high-dimensional classification problem; Fan, Peng and Huang (2005) and Huang \etal(2005) studied semiparametric inference where $\nu\rai$; and Hu and He (2006) proposed an enhanced quantile normalization based on high-dimensional singular value decomposition to reduce information loss in gene expression profiles.

Korosok and Ma (2005) treated the problem of uniform, simultaneous estimation of a large number of marginal distributions, showing that if $\log\nu=o(n)$, and if certain other conditions hold, then $\max_{1\leq i\leq\nu}\,\|\hF_i-F_i\|_\infty\ra0$, where $\hF_i$ is an estimator of the $i$th marginal distribution~$F_i$.  As a corollary they proved that a P-value $\hP_i$ of $\hF_i$ converges uniformly in $i$ to its counterpart $P_i$ for $F_i$, provided $\log \nu=o(n\half)$:
$$
\max_{1 \leq i\leq \nu}\,\|\hP_i-P_i\|_\infty \to 0\,.
\eqno{(1.1)}
$$
These results are important advances in the literature of simultaneous testing, where P-values are popularly assumed to be known. For examples in the latter setting, see Benjamini and Yekutieli (2001), Dudoit, Shaffer and Boldrick (2003), 
Donoho and Jin (2004), Efron (2004), Genovese and Wasserman (2001), Storey, Taylor and Siegmund (2004), Lehmann and Romano (2005), Lehmann, Romano and Shaffer (2005) where many new ideas have been introduced to control different aspects of simultaneous hypothesis testing and false discovery rate~(FDR).

In many practical settings the assumption that P-values are
calculated without error is unrealistic, but it is unclear how good
the approximation must be in order for simultaneous inference to be
feasible.  Simple consistency, as evidenced by (1.1), is not enough;
the level of accuracy required must increase with~$\nu$. More
precisely, letting $\alpha_\nu $ be the significant level, which
tends to zero as $\nu \to \infty$, the required accuracy is then
$$
\max_{1 \leq i\leq \nu}\,\|\hP_i - P_i \|_\infty = o(\alpha_\nu).
\eqno{(1.2)}
$$
In this paper we provide a concise solution to this problem.

For example, we show that in the case of simultaneous $t$-tests,
calibrated by reference to Normal or Student's~$t$ distributions, a
necessary and sufficient condition for overall level accuracy to be
asymptotically correct is (a)~$\log\nu=o(n^{1/3})$.  This is true
even if the sampling distribution is highly skewed or heavy tailed.
On the other hand, if bootstrap methods are used for estimating
P-values then the asymptotic level of the simultaneous test is
accurate as long as (b)~$\log \nu = o(n\half)$. These results make
clear the advantages offered by bootstrap calibration.  We shall
explore them numerically as well as theoretically.  Result (a) needs
only bounded third moments of the sampling distribution, although
our proof of (b) needs more restrictions.

Take the case of family-wise error rate as an example. If the
overall error rate is controlled at $p$, then $k_n$ hypotheses with
the smallest P-values are rejected, where
$$
k_n=\max\{i:P_{i} \leq  p/\nu,i=1,\ldots,\nu\}\eqno{(1.3)}
$$
and $P_i$ denotes the significance level of the $i$th test. The
Benjamini and Hochberg's (1995) approach to control false discovery
rate (FDR) at $p$ is to select
$$
k_n=\max\{i:P_{(i)} \leq  i p/\nu,i=1,\ldots,\nu\}, \eqno{(1.4)}
$$
where $\{P_{(i)}\}$ are the ordered values of $\{P_i\}$. If the
distributions from which the $P_i$'s are computed need to be
estimated then, in view of (1.3) or (1.4), the error of the
estimators $\hP_i$ should equal $o(\nu\mo)$ in order to sort
correctly $\{P_i\}$, and the approximation (1.1) requires
significant refinement.  Indeed, it corresponds to (1.2) with
$\alpha_\nu = p/\nu$. This is a very stringent requirement
($\alpha_\nu = 10^{-5}$, if $p=0.1$ and $\nu = 10^4$) and the
accuracy is hard to
achieve for many practical situations. 

Sometimes, the classical approach provides an attractive alternative
to select significant hypotheses (genes).  For example, in their
analysis of gene expression data, Fan \etal(2004) take
$\a_\nu=0.001$ and find the significant set of genes,
$$
{\cal S}=\{i:P_i\leq\a_\nu,i=1,\ldots,\nu\}\,, \eqno{(1.5)}
$$
for $\nu=15,\!000$ simultaneous tests.  Here $\a_\nu$ is an order of
magnitude larger than $\nu^{-1}$, and the approximation errors when
estimating P-values need only be $o(\a_\nu)$ when computing (1.5),
rather than $o(\nu\mo)$ in the family-wise error rate problem. In
this case, the requirement is much less stringent and the number of
simultaneous tests $\nu$ for (1.2) to hold should be an order of
magnitude larger than the case with $\alpha_\nu = \nu^{-1}$.
However, even if we take $\a_\nu=C\,\nu^{-a}$, where $C>0$ and $a
\in (0,1]$ are constants, from some points of view the problem does
not become appreciably simpler.  For instance, in the case of
calibration by the Normal or Student's~$t$ distributions, condition
(a) two paragraphs above is still necessary and sufficient.

In this example the number of elements of ${\cal S}$, denoted by
$k_n'$, equals the number of genes discovered at the significance
level $\a_\nu$; and $\nu \a_\nu$, an upper bound to the expected
number of falsely discovered genes, is approximately the same as the
expected number of falsely discovered genes when most null
hypotheses are true. Hence, the FDR in this example is estimated as
$\hat{p}=\nu\a_\nu/k_n'$.  Surprisingly, this classical approach is
approximately the same as the Benjamini-Hochberg method.  Using
$\hat{p}$ as the FDR to be controlled, the Benjamini-Hochberg method
picks up $\hat{k}_n =\max\{i:P_{(i)} \leq
i\hat{p}/\nu,i=1,\ldots,\nu\}$ most significant P-values, whereas
the classical approach (1.5) selects $k_n'$ most significant
P-values. The classical approach is more conserve since it is shown
by Fan \etal (2005) that $k_n' \leq \hat{k}_n$. However, in our
simulation experiments (reported in a previous version), $\hat{k}_n
\approx k_n$. In other words, the classical approach selects nearly
the same set of significant hypotheses as that by the
Benjamini-Hochberg method.

The results stated above are for the case of independent tests, but
they also apply, in the sense of sufficiency for asymptotically
conservative tests, under the assumption of positive regression
dependency; see Benjamini and Yekuteli (2001, p.~1170) for
discussion.  We discuss too the case where the errors in the data
come from a linear process.  There, theoretical formulae for
significance levels can be based on Poisson cluster-process
arguments.  To treat simultaneous hypothesis testing under
completely general types of dependence we use a Bonferroni argument.
This permits our results on the accuracy of bootstrap and Student's
$t$ approximations to be extended to a wide range of settings.  In
particular they apply in the context of generalized family-wise
error rate.

Practical implications of our work include the following. (i)~When
calibrating multiple hypothesis tests using Student's $t$ statistic,
for example with a view to controlling family-wise error rate,
impressive level-accuracy can be obtained via Student's $t$
approximation.  Only a mild moment condition is needed.
(ii)~Nevertheless, accuracy is noticeably improved through using the
bootstrap.  (iii)~These results also apply in the case of
generalized family-wise error rate. (iv)~Owing to the limitation of
on the number of simultaneous hypotheses that can be accurately
tested for a given $n$ and $\alpha_\nu$, other methods of
pre-screening are needed when there are excessively many hypotheses
to be tested.

In a sequence of papers, Finner and Roters (1998, 1999, 2000, 2001, 2002) developed theoretical properties of $n$ simultaneous hypothesis tests as $n$ increases.  However, their work differs from ours in a major respect, through their assumption that the exact significance level can be tuned to a known value in the continuum.  In such cases there is no theoretical limit to how large $n\nu$ can be.  By way of comparison, in the setting of the genetic problems that motivates our work, level inaccuracies limit the effective size of $\nu$; we shall delineate this limitation using both theoretical and numerical arguments.

The paper is organized as follows.  In section~2 we formulate the accuracy problem for simultaneous tests. There, we also outline statistical models and testing procedures.  Our main results are presented in section~3, where we answer the question of how many hypotheses can be tested simultaneously without the overall significance level being seriously in error.  The theoretical definition of the latter property is that the overall significance level should converge to the nominal one as the number of tests increases, if the samples on which the tests are based are independent; and that the limiting level should not exceed the nominal one if the independence condition is violated and a Bonferroni bound is used.

Section~4 outlines the idea of marginal aggregation when the number of hypotheses is ultra-large.  Numerical investigations among various calibration methods are presented in section~5.  Technical proofs of results in section~3 are relegated to section~6.

\bs

\cl{\bf 2.  MODEL AND METHODS FOR TESTING}

\ni{\sl 2.1.  Basic model and methodology.}
The simplest model is that where we observe random variables
$$
Y\ij=\mu_i+\ep\ij\,,\quad1\leq i<\infty\,,
\quad1\leq j\leq n\,,\eqno(2.1)
$$
with the index $i$ denoting the $i$th gene, $j$ indicating the $j$th
array, and the constant $\mu_i$ representing the mean effect for the
$i$th gene.  We shall assume that:
$$
\eqalign{ &\hbox to4.4in{for each $i$, $\ep_{i1},\ldots,\ep_{in}$
are independent and identically distributed random}\cr
\noalign{\vskip-7pt} &\hbox{variables with zero expected
value.}\cr}\eqno(2.2)
$$
The results given below are readily extended to the case where $n=n_i$ depends on $i$, but taking $n$ fixed simplifies our discussion.

Let $T_i=n\half\,\bY_i/S_i$, where
$$
\bY_i=\oon\,\sumjon Y\ij\,,\quad
S_i^2=\oon\,\sumjon(Y\ij-\bY_i)^2\,.
$$
For a given value of $i$ we wish to test the null hypothesis $H\zi$
that $\mu_i=0$, against the alternative hypothesis $H\oi$ that
$\mu_i\neq0$, for $1\leq i\leq\nu$ say. We first study this
classical testing problem of controlling the probability of making
at least one false discovery, which requires calculating P-values
with accuracy $o(\nu^{-1})$, the same as that needed in (1.3). We
then extend our results to control the relaxed FDR in (1.5), which
is less stringent.

A standard test is to reject $H\zi$ if $|T_i|>t_\a$.  Here, $t_\a$
denotes the solution of either of equations
$$
P(|Z|>t_\a)=1-(1-\a)^{1/\nu}\,,\quad
P\{|T(n-1)|>t_\a\}=1-(1-\a)^{1/\nu}\,,\eqno(2.3)
$$
where $N$ and $T(k)$ have respectively the standard Normal
distribution and Student's $t$ distribution with $k$ degrees of
freedom. Note that (2.3) serves only to give a definition of $t_\a$
that is commonly used in practice; it does not amount to an
assumption about the sampling distribution of the data. Indeed,
$\alpha_\nu = 1 - (1 - \alpha)^{1/\nu}$ and $\alpha$ is a one-to-one
map. The core of the argument in this paper is that the accuracy of
the distributional approximations implicit in (2.3) are based on a
delicate relationship between $n$ and $\nu$, which is central to the
question of how many simultaneous tests are possible.

\ni{\sl 2.2.  Significance levels for simultaneous tests.}
If $H\zi$ is true then the significance level of the test restricted
to gene $i$, is given by
$$
\p_i=P\zi(|T_i|>t_\a)\,,\eqno(2.4)
$$
where $P\zi$ denotes probability calculated under~$H\zi$. For the
classical approach (1.3), which is nearly the same as the
Benjamini-Hochberg method, we ask how large $\nu$ can be so that
$$
    \max_{1 \leq i \leq \nu} | p_i - \alpha_\nu| = o(\alpha_\nu)?
    \eqno{(2.5)}
$$
The answer depends on the rate at which $\alpha_\nu$ approaches to
zero.  In particular, if $\alpha_\nu = \beta / \nu$, (2.5) entails
that
$$
\max_{1\leq i\leq\nu}\,\p_i=o(1)\quad\hbox{and}\quad
\sumionu\p_i=\be+o(1)\,,\eqno(2.6)
$$
for some $0 < \beta < \infty$?  Result (2.6) implies that the
significance level of the simultaneous test, described in
section~2.1, is
$$
\eqalignno{ \a(\nu) &\equiv P\Big(\hbox{$H\zi$ rejected for at least
one $i$ in the range $1\leq i\leq\nu$}\Big)&(2.7)\cr &
 \leq\sumionu\p_i =\be+o(1) \,.\qquad&(2.8)\cr}
$$
If, in addition to (2.2), we assume that
$$
\hbox to4.4in{the sets of variables $\{\ep\ij$, $1\leq j\leq n\}$
 are independent for different $i$,}\eqno(2.9)
$$
then
$$
\eqalignno{ \a(\nu) &=1-\prod_{i=1}^\nu\,(1-\p_i)
=1-\exp\bigg(-\sumionu\p_i\bigg)+O\bigg(\sumionu\p_i^2\bigg)
.\qquad&(2.10)\cr}
$$
Consequently, (2.6) and (2.10) imply the following property:
$$
\hbox{if (2.9) holds then\quad
$\a(\nu)=1-e^{-\be}+o(1)\,,$}\eqno(2.11)
$$
where $\a(\nu)$ is as defined at~(2.7). The ``$o(1)$'' terms in
(2.8) and (2.11) are quantities which converge to zero as $\nu\rai$.
Result (2.11) also holds, with the identity
$\a(\nu)=1-e^{-\be}+o(1)$ replaced by $\a(\nu)\leq1-e^{-\be}+o(1)$,
if (2.9) is replaced by the positive regression dependency
assumption (Benjamini and Yekuteli, 2001 p.~1170).

In practice we would take $\be=-\log(1-\a)$, if we were prepared to
assume (2.9) (or positive regression dependency of the test
statistics) and wished to construct a simultaneous test with level
close to~$\a$ [or, respectively, asymptotically not exceeding~$\a$];
and take $\be=\a$, if we were using Bonferroni's bound to construct
a conservative simultaneous test with the same approximate level.
Further discussion of the dependent-data case is given in
section~2.3.

\ni{\sl 2.3.  Generalized family-wise error rate.} The results above
can be generalized by extending the definition at (2.7) to,
$$
\eqalignno{
\a_k(\nu) &\equiv P\Big(\hbox{$H\zi$ rejected for at least
$k$ values of $i$ in the range $1\leq i\leq\nu$}\Big)\cr
&
 \leq{1\over k}\,\sumionu\p_i  =k\mo\be+o(1) \,,\qquad&(2.12)\cr}
$$
where (2.12) follows from (2.6).  We also have the following
analogue of (2.11): Assuming~(2.6):
$$
\hbox{if (2.9) holds then}\quad
\a_k(\nu)=1-\sum_{j=1}^k\,{\be^j\over
j!}\;e^{-\be}+o(1)\,.\eqno(2.13)
$$
Under the positive regression dependency assumption, the equality in
(2.13) would be replaced by~$\leq$.

Insight into properties of generalized family-wise error-rate in
cases of dependency can be gained by considering settings where the
processes $\cP_j=\{\ep\ij,i\geq1\}$, for $j\geq1$, are independent
and identically distributed moving averages, each with the
distribution of $\cP=\{\ep_i,i\geq1\}$, where
$\ep_i=\sumk\,\th_k\,\de_{i+k}$, the variables $\de_i$ are
independent and identically distributed with zero mean, and the
weights $\th_k$ may depend on $n$ but satisfy $\sumk\th_k^2\leq C$
for a fixed constant~$C$.  Here, under conditions on the
distribution of $\de$, it can be shown that if (2.6) holds then the
limiting distribution of $N$, defined to equal the number of indices
$i$ for which $|T_i|>t_\a$, is that of $\sum_{1\leq\ell\leq
M}\,K_\ell$, where $M$ has a Poisson distribution with mean
$\la=\be/E(K)$, $K_1,K_2,\ldots$ are identically distributed as the
nonnegative, integer-valued random variable $K$, and
$M,K_1,K_2,\ldots$ are independent.  Thus, the distribution of the
number of counts is based on a Poisson cluster process, with $K$
denoting cluster size.

In this case, without the independence assumption~(2.9),
$$
\a_k(\nu)\ra\sum_{m=0}^\infty\,{1\over m!}\;e^{-\la}\,
P\bigg(\sum_{\ell=1}^m\,K_\ell\geq k\bigg)\leq k\mo\,\be\,,
$$
which property generalizes both (2.12) and~(2.13).  This result is
not necessarily beneficial in practice, however, owing to the
difficulty of estimating the distribution of~$K$.

\ni{\sl 2.4.   Methods for calibration.}
For calibration against Normal or Student's $t$ distributions we
take the critical point $t_\a$ to be the solution of the respective
equations~(2.3).  Below we consider bootstrap calibration; Edgeworth correction (see e.g.~Hall, 1990) could also be used.

Let $Y_{i1}\dag,\ldots,Y_{in}\dag$ denote a bootstrap resample drawn
by sampling randomly, with replacement, from
$\cY_i=\{Y_{i1},\ldots,Y_{in}\}$.  Put $Y\ij\as=Y\ij\dag-\bY_i$ and
$T_i\as=n\half\bY_i\as/S_i\as$, where $\bY_i\as=n\mo\,\sumj Y\ij\as$
and $(S_i\as)^2=n\mo\,\sumj(Y\ij\as-\bY_i\as)^2$.  Write $z_\a$ for
the conventional Normal critical point for $\nu$ simultaneous tests.
That is, $z_\a$ solves the equation $P(|Z|>z_\a)=1-(1-\a)^{1/\nu}$.
(We could also use the Student's $t$ point.)  Define $\xi=\hf_i(\a)$
to be the solution of the equation
$$
P(|T_i\as|>z_\xi\mi\cY_i)=1-(1-\a)^{1/\nu}\,.
$$
Our bootstrap critical point is $\hatt\ia=z_{\hf_i(\a)}$; we reject $H\zi$ if and only if~$|T_i|>\hatt\ia$.  The definition of $\p_i$ at (2.4) should here be replaced by,
$$
p_i=P\zi(|T_i|>\hatt\ia)\,.\eqno(2.14)
$$
With this new definition, (2.11), (2.12) and (2.13) continue to be a
consequences of~(2.6).

\bs

\cl{\bf 3.  THEORETICAL RESULTS}

\ni{\sl 3.1.  Asymptotic results.}
Define $\ka_{i3}$ to be the third cumulant, or equivalently the skewness, of the distribution of $\ep_i'=\ep_{i1}/(E\ep_{i1}^2)\half$.

\proclaim Theorem~3.1.  Assume that
$$
\max_{1\leq i\leq\nu}\,E|\ep_i'|^3=O(1)\eqno(3.1)
$$
as $\nu\rai$, and suppose too that $\nu=\nu(n)\rai$ in such a manner
that $(\log\nu)/n^{1/3}\ab\ra\ga$, where $0\leq\ga<\infty$.  Define
$t_\a$ by either of the formulae at $(2.3)$, and $\p_i$ by~$(2.4)$.
Then $(2.6)$ holds with
$$
\be=\be(\nu)\equiv-{\log(1-\a)\over\nu}\;
\sumionu\cosh\big(\otd\,\ga^3\ka_{i3}\big)\, , \eqno(3.2)
$$
where $\cosh(x)=(e^x + e^{-x})/2$.

The value of $\be(\nu)$, defined at (3.2), is bounded by
$|\log(1-\a)|\,\cosh(\ga^3\,B)$, uniformly in $\nu$, where $B=\sup_i|\ka_{i3}|$.

\proclaim Corollary~3.2.  Assume the conditions of Theorem~3.1.  If
$\ga=0$, i.e.~if $\log\nu=o(n^{1/3})$, then $(2.6)$ holds with
$\be=-\log(1-\a)$; and if $\ga>0$ then $(2.6)$ holds with
$\be=-\log(1-\a)$ if and only if
$\nu\mo\,\sum_{i\leq\nu}\,|\ka_{i3}|\ra0$, i.e.~if and only if the
limit of the average absolute values of the skewnesses of the
distributions of $\ep_{11},\ldots,\ep_{\nu1}$ equals zero.

\ni Since Corollary~3.2 implies (2.6) then it also entails (2.8) and
(2.11)--(2.13).

\proclaim Theorem~3.3.  Strengthen $(3.1)$ to the assumption that
for a constant $C>0$, $P(|\ep_i'|\leq C)=1$, and suppose too that
$\nu=\nu(n)\rai$ in such a manner that $\log\nu=o(n\half)$.  Define
$\hatt\ia=z_{\hf_i(\a)}$, as in section~2.4, and define $\p_i$
by~$(2.14)$.  Then $(2.6)$ holds with $\be=-\log(1-\a)$.

\ni{\sl 3.2. Applications to controlling error rate.}
Define $t_\a$ and $\hatt\ia$ by (2.3) and as in section~2.4, respectively.  In the proof of Theorem~3.1 it is shown that, with $\be=-\log(1-\a)$ and using conventional calibration, $P\zi(|T_i|>t_\a)=\be\,\nu\mo+o(\nu^{-1})$, uniformly in $i$ under the null hypotheses, provided $\log\nu=o(n^{1/3})$; and that, when employing bootstrap calibration, $P\zi(|T_i|>\hatt\ia)=\be\,\nu\mo+o(\nu\mo)$, again uniformly in $i$, if $\log\nu=o(n^{1/2})$.  These results substantially improve a uniform convergence property of Kosorok and Ma (2005), at the expense of more restrictions on~$\nu$.

When the P-values in (1.3) need to be estimated, the estimation
errors should be of order $o(\nu^{-1})$, where $\nu$ diverges
with~$n$. On the other hand, when P-values in (1.5) are estimated,
the precision can be of order $o(\a_\nu)$, where for definiteness we
shall take $\a_\nu=C\,\nu^{-a}$ with $C>0$ and $a \in (0,1]$.  In
the latter case the results in Theorems 3.1 and 3.3 continue to
apply; there is no relaxation, despite the potential simplicity of
the problem.

To appreciate why, note that the tail probability of the standard
Normal distribution satisfies $P(|Z|\geq
x)\sim\exp(-x^2/2)/(\sqrt{2\pi}\,x)$. Suppose that the large
deviation result holds up to the point $x=x_n$, which should be of
order $o(n^{1/6})$ for Student's $t$ calibration, and $o(n^{1/4})$
for bootstrap calibration. Setting $P(|Z|\geq x)$ equal to $\a_\nu$
yields $\log(1/\a_\nu)=\thf\,x_n^2+\log x_n+\ldots$, that is,
$$
a \,\log\nu=\thf\,x_n^2+\log x_n+\log(\sqrt{2\pi}\,C) +\hbox{smaller
order terms}\,.
$$
Regardless of the values of $C>0$ and $a\in (0,1]$, this relation
implies that the condition $x_n=o(n^{1/6})$ is equivalent to
$\log\nu=o(n^{1/3})$ and $x_n=o(n^{1/4})$ entails
$\log\nu=o(n^{1/2})$, although taking $a$ close to~0 will
numerically improve approximations, in both theory and practice.

\bs

\cl{\bf 4.  MARGINAL AGGREGATION}

We have shown that with bootstrap calibration, Student's
$t$-statistics can test simultaneously a number of hypotheses of
order $\exp\{o(n\half)\}$. Although this value may be conservative,
it may still not be large enough for some applications to microarray
and tiling arrays where the number of simultaneous tests can be even
larger. Similarly, whether we consult a $t$-table or a Normal table
depends very much on mathematical assumptions.  For example, suppose
an observed value of a $t$ statistic is 6.7. Its corresponding
two-tail P-value, for $n=6$ arrays, is 0.112\% when looking up
$t$-tables with five degrees of freedom, and $2.084 \times 10^{-11}$
when consulting Normal tables.   Yet, both distributions can be
regarded as the approximate ones.

To overcome these problems, Reiner \etal(2003) and Fan \etal(2004)
introduce marginal aggregation methods. The basic assumption is that
the null distributions of test statistics $T_i$ are the same,
denoted by $F$. With this assumption (which also implicitly assumed
when the same distribution table is looked up), we can use the
empirical distribution of $\{T_i,i=1, \ldots,\nu\}$,~i.e.
$$
\hF_\nu (x)=\nu^{-1}\,\sum_{i=1}^\nu\,I(T_i \leq x)\,,\eqno{(4.1)}
$$
as an estimator of $F$.  This turns the ``curse-of-dimensionality''
into a ``blessing-of-dimensionality''. In this case, even if $n$ is
fixed, the distribution of $F$ can still be consistently estimated
when $\nu \to \infty$ and only a small fraction of alternative
hypotheses are true.

Reiner \etal(2003) and Fan \etal(2004) carry this idea one step further.  They aggregate the estimated distributions for each $T_i$, based on a resampling technique (more precisely, a permutation method). For example we can use, as an estimator of $F(x)$, the average of bootstrap estimators,
$$
\hF_\nu^*(x)=\nu^{-1}\,\sum_{i=1}^\nu\,
P(T_i\as\leq x\mi{\cal Y}_i )\,.\eqno{(4.2)}
$$

The theorem below describes properties of $\hF_\nu$ and
$\hF_\nu\as$, defined at (4.1) and (4.2), respectively.  Interpret
$x_n$, below, as a left-tail critical value; similar results also
hold in the upper tail.

\proclaim Theorem~4.1.  Let $\nu_1$ be the number of nonzero $\mu_i$, i.e.~the number of elements of the set $\{i: H_{1i}$ is true$\}$, and put $F_n(x)=E\{P(T_i^* < x)\}$. Then,
$$
\eqalignno{ \hF_\nu (x_n)&=F(x_n)+O_p\left[(\nu_1/\nu)+
\sqrt{F(x_n)/\nu}\;\big\{1+a_n\,\nu
+a_n\,\nu_1\,F(x_n)\mhf\big\}\half \right]\,,\cr
\hF_\nu\as(x_n)&=F_n(x_n)+O_p\left \{ \sqrt{F(x_n)/\nu} \right
\}\,,\cr}
$$
provided that $|r_{ij}| \leq a_n$, with $r_{ij}$ denoting the correlation coefficient between $I(T_i\leq x_n)$ and $I(T_j\leq x_n)$.

The term $O(\nu_1/\nu)$ in (4.3) reflects the bias of the estimate.  It can be reduced by using the sieve idea of Fan \etal(2005).

\bs

\cl{\bf 5.  NUMERICAL PROPERTIES}

In our simulation study we construct models that reflect aspects of gene expression data. To this end, we divide genes into three groups. Within each group, genes share one unobserved common factor with different factor loadings.  In addition, there is an unobserved common factor among all the genes across the three groups. For simplicity of presentation, we assume that $\nu$ is a multiple of three.  We denote by $Z_{ij}$ a sequence of independent N$(0, 1)$ random variables, and $\chi_{ij}$ a sequence of independent random variables of the same distribution as that of $(\chi_{m}^2 - m)/\sqrt{2m}$.  Note that $\chi_{ij}$ has mean 0, variance~1 and skewness $\sqrt{8/m}$.  In our simulation study we set $m=6$.

With given factor loading coefficients $a_i$ and $b_i$, the error
$\ep_{ij}$ in (2.1) is defined as
$$
\ep_{ij}={Z_{ij}+a_{i1}\,\chi_{j1}+\a_{i2}\,\chi_{j2}+
a_{i3}\,\chi_{j3}+b_i\,\chi_{j4}\over
(1+a_{i1}^2+a_{i2}^2+a_{i3}^2+b_i^2)\half}\,,\quad
i=1,\ldots,\nu\,,\quad j=1,\ldots,n\,,
$$
where $a_{ij}=0$ except that $a_{i1}=a_i$ for $i=1,\ldots,\otd\,\nu$, $a_{i2}=a_i$ for $i=\otd\,\nu+1, \ldots,\tthd\,\nu$, and $a_{i3}= a_i$ with $i=\tthd\,\nu+1,\ldots,\nu$.  Note that $E\ep_{ij}= 0$ and $\var(\ep_{ij})=1$, and that the within-group correlation is in general stronger than the between-group correlation, since the former shares one extra common factor.  We consider two specific choices of factor loadings: Case~I, where the factor loadings are taken to be $a_j=0.25$ and $b_j=0.1$ for all $j$ (thus, the $\ep_{ij}$'s have the same marginal distribution, although they are correlated); and Case~II, for which the factor loadings $a_i$ and $b_i$ are generated independently from, respectively, $U(0,0.4)$ and $U(0,0.2)$.

The ``true gene expression'' levels $\mu_i$ are taken from a realization of the mixture of a point mass at~0 and a double-exponential distribution: $c\,\delta_{0}+\thf\,(1-c)\,\exp(-|x|)$, where $c \in (0, 1)$ is a constant.  With the noise and the expression level given above, $Y_{ij}$ generated from (2.1) represents, for each fixed $j$, the observed log-ratios between the two-channel outputs of a c-DNA microarray. Note that $|\mu_j| \geq \log2$ means that the true expression ratio exceeds~2.  The probability, or the empirical fraction, of this event equals~$\thf\,(1-c)$.

For each $\a_\nu$ we compute the P-value according to the Normal
approximation, $t$-approximation, the bootstrap method and the
aggregated bootstrap~(4.2).  This results in $\nu$ estimated
P-values $\hP_j$ for each method and each simulation. Let $N_1$
denote the number of P-values that are no larger than $\a_\nu$;
see~(1.5). Then, $N_1/\nu$ is the empirical fraction of null
hypotheses that are rejected.  When $c=1$, $N_1/(\nu\a_\nu)-1$
reflects the accuracy with which we approximate P-values.  Its root
mean square error (RMSE), $\{E(N_1/(\nu\a_\nu)-1)^2\}\half$, will be
reported, where the expectations are approximated by averages across
simulations. We exclude the marginal aggregation method~(4.1), since
in our simulations it always gave $N_1/\nu=\a_\nu$.

We take $\nu=600$ (small), $\nu=1,800$ (moderate) and $\nu=6,000$ (typical) for microarray applications (after preprocessing, which filters out many low quality measurements on certain genes) and $\a_\nu=1.5 \nu^{-2/3}$, resulting in $\a_\nu=0.02$, $0.01$ and $0.005$, respectively.  The sample size $n$ is taken to be 6 (typical number of microarrays), 20 (moderate) and 50 (large).  The number of replications in simulations is $600,\!000/\nu$.  For the bootstrap calibration method and the aggregated method (4.2), we replicate bootstrap samples 2,000, 4,000 and 9,000 times, respectively, for $\a_\nu=0.02, 0.01$ and 0.005.

Tables~1 and~2 report the accuracy of estimated P-values when $c=1$.
It can be seen that the Normal approximations are too inaccurate to
be useful. Therefore we shall exclude the Normal method in the
discussion below.  For $n=20$ and 50, the bootstrap method provides
better approximations than Student's $t$-method. This indicates that
the bootstrap can test more hypotheses simultaneously, which is in
accord with our asymptotic theory. Overall the bootstrap method is
also slightly better then the aggregated bootstrap (4.2), although
the two methods are effectively comparable. However, with the small
sample size $n=6$,  Student's $t$-method is relatively the best,
although the approximations are poor in general. This is
understandable, as the noise distribution is not Normal. With such a
small sample size, the two bootstrap-based methods, in particular
the aggregated bootstrap method (4.2), suffer more from random
fluctuation in the original samples.

We now illustrate the four methods by using a microarray data set.
The data were analyzed by Fan \etal(2004), where the biological aim
was to examine the impact of the stimulation by MIF, a growth
factor, on the expressions of genes in neuroblastoma cells. Six
arrays of cDNA microarray data were collected, consisting of
relative expression profiles of 19,968 genes in the MIF stimulated
neuroblastoma cells (treatment) and those without stimulation
(control). After preprocessing that filtered low quality expression
profiles, 15,266 genes remained.  Within-array normalization,
discussed by Fan \etal(2004), was applied to remove the intensity
effect and block effect. Among 15,266 gene expression profiles, for
simplicity of illustration we focussed only on 7,583 genes that do
not have any missing values.  In our notation, $\nu= 7,\!583$ and
$n=6$. Table~3 summarizes the results at different levels of
significance.

Different
methods for estimating P-values yield very different results. In particular, the
distribution of P-values computed by looking up the normal table is
stochastically much larger than that based on the $t$-table, which
in turn is stochastically much larger than that based on the
bootstrap method.  The results are very different. As noted before,
the normal approximation is usually very poor and grossly inflates
the number of significant genes.  The same remark applies to the
$t$-approximation with moderate or large numbers of degrees of freedom. The
most accurate approximation is given by bootstrap methods.

\bs

\cl{\bf 6.  PROOFS OF RESULTS IN SECTION~3}

For the sake of brevity we shall derive only Theorems~3.1 and~3.3.  Let $C_1>0$.  Given a random variable $X$ with $E(X)=0$, consider the condition:
$$
E(X)=0\,,\quad E\big(X^2\big)=1\,,\quad E\big(X^4\big)\leq C_1\,.\eqno(6.1)
$$
The following result follows from Theorem~1.2 of Wang (2005), after
transforming the distribution of $T$ to that of $(\sumi X_i)/(\sumi X_i^2)$.

\proclaim Theorem~6.1.  Let $X,X_1,X_2,\ldots$ denote independent and identically distributed random variables such that $(6.1)$ holds.  Write $T=T(n)$ for Student's $t$ statistic computed from the sample $X_1,\ldots,X_n$, with (for the sake of definiteness) divisor $n$ rather than $n-1$ used for the variance.  Put $\pi_3=-\otd\,\ka_3$, where $\ka_3$ denotes the skewness of the distribution of $X/(\var X)\half$.  Then,
$$
{P(T>x)\over1-\Phi(x)}
=\exp\big(\pi_3\,x^3\,n\mhf\big)\,
\bigg\{1+\th\;{(1+x)^2\over n\half}\bigg\}\,,
\eqno(6.2)
$$
where $\th=\th(x,n)$ satisfies $|\th(n,x)|\leq C_2$ uniformly in $0\leq x\leq C_3\,n^{-1/4}$ and $n\geq1$, and $C_2,C_3>0$ depend only on $C_1$.

Theorem~3.1 in the case of Normal calibration follows directly from
Theorem~6.1. The case of Student's $t$ calibration can be treated
similarly.

To derive Theorem~3.3, note that each $\var(\ep_i')=1$.  To check that, with probability at least $p_n\equiv1-\exp(-d_1\,n\half)$ for a constant $d_1>0$, the conditions of Theorem~6.1 hold for the bootstrap distribution of the statistic $T_i^*$, for each $1\leq i\leq\nu$, it suffices to show that there exist constants $0<C_4<C_5\half$ such that, with probability at least $p_n$, the following condition holds for $1\leq i\leq\nu$:
$$
C_4\leq\oon\,\sumjon(Y_{ij}-\bY_i)^2\,,\quad
\oon\,\sumjon(Y_{ij}-\bY_i)^4\leq C_5\,.\eqno(6.3)
$$
This can be done using Bernstein's inequality and the assumption
that, for each $i$, $P(|\ep_i'|\leq C)=1$, and can also be shown by
the uniform convergence result of the empirical process of Korosok
and Ma (2005).

Let $\cE_n$ denote the event that (6.3) holds for each~$1\leq i\leq\nu$.  When $\cE_n$ prevails, we may apply Theorem~6.1 to the distribution of $T_i\as$ conditional on $\cY_i$, obtaining:
$$
P\big(T_i\as>x\bigmi\cY_i\big) =\{1-\Phi(x)\}\,
\exp\big(-\otd\,\hka_{i3}\,n\mhf\,x^3\big)\,
\bigg\{1+\hth_i\;{(1+x)^2\over n\half}\bigg\}\,,\eqno(6.4)
$$
where $\hka_{i3}$ is the empirical version of $\ka_{i3}$, computed from $\cY_i$, and, on an event of which the probability equals $1-O\{\exp(-d_2\,n\half)\}$, $|\hth_i|\leq D_1$ uniformly in $i$ and in $0\leq x\leq x_n$.  (Here and below, $x_n$ will denote any sequence diverging to infinity but satisfying~$x_n=o(n^{1/4})$, and $D_1,D_2,\ldots$ and $d_1,d_2,\ldots$ will denote constants.)  It follows directly from Theorem~6.1 that
$$
P\zi(T_i>x) =\{1-\Phi(x)\}\,
\exp\big(-\otd\,\ka_{i3}\,n\mhf\,x^3\big)\,
\bigg\{1+\th\;{(1+x)^2\over n\half}\bigg\}\,,\eqno(6.5)
$$
where $|\th_i|\leq D_1$ uniformly in $i$ and in $0\leq x\leq x_n$.

Result (6.5), and its analogue for the left-hand tail of the
distribution of $T_i$, allow us to express $t\ia$, the solution of
the equation $P\zi(|T_i|>t\ia)=1-(1-\a)^{1/\nu}$, as a Taylor
expansion:
$$
\big|t\ia-z_\a-c\,\ka_{i3}\,n\mhf\,z_\a^2\big|
\leq D_2\,\big(n\mo\,z_\a^4+n\mhf\,z_\a\big)\,,
$$
uniformly in $i$, where $c$ is a constant and $z_\a$ is the solution
of $P(|Z|>z_\a)=1-(1-\a)^{1/\nu}$.  Note that if $z_\a$ solves this
equation then $z_\a\sim(2\,\log\nu)\half$, and so, since
$\log\nu=o(n\half)$, then $z_\a=o(n^{1/4})$.  Therefore, without
loss of generality, $0\leq z_\a\leq x_n$.  Likewise we may assume
below that $0\leq t_{i\a}\leq x_n$, and $0\leq\hatt_{i\a}\leq x_n$
with probability $1-O\{\exp(-d_2\,n\half)\}$.

Also, from (6.4) we can see that on an event of which the probability equals $1-O\{\exp(-d_2\,n\half)\}$,
$$
\big|\hatt\ia-z_\a-c\,\hka_{i3}\,n\mhf\,z_\a^2\big|
\leq D_3\,\big(n\mo\,z_\a^4+n\mhf\,z_\a\big)\,.
$$
However, on an event with probability
$1-O\{\exp(-d_3\,n\half)\}$, $|\hka_{i3}-\ka_{i3}|\leq
D_4\,n^{-1/4}$, and therefore, on an event with probability
$1-O\{\exp(-d_4\,n\half)\}$,
$$
\big|\hatt\ia-z_\a-c\,\ka_{i3}\,n\mhf\,z_\a^2\big|
\leq D_5\,\big(n\mo\,z_\a^4+n\mhf\,z_\a
+n^{-3/4}\,z_\a^2\big)\,.
$$

It follows from the above results that $P\zi(|T_i|>\hatt\ia)$ lies
between the respective values of
$$
P\zi(|T_i|>t\ia\pm\de) \mp
D_6\,\exp\big(-d_4\,n\half\big)\,,\eqno(6.6)
$$
where
$$
\de=D_5\,\big(n\mo\,z_\a^4+n\mhf\,z_\a
+n^{-3/4}\,z_\a^2\big)\,.
$$
Using (6.5), and its analogue for the left-hand tail, to expand the
probability in (6.6), we deduce that
$$
P\zi(|T_i|>t\ia\pm\de)=P\zi(|T_i|>t\ia)\,\{1+o(1)\}\,,
$$
uniformly in $i$.  More simply,
$\exp(-d_4\,n\half)=o\{P\zi(|T_i|>t\ia)\}$, using the fact that $z_\a=o(n^{1/4})$ and $\exp(-D_7\,z_\a^2)=o\{P\zi(|T_i|>t\ia)\}$ for sufficiently large
$D_7>0$.  Hence,
$$
P\zi(|T_i|>\hatt\ia)=P\zi(|T_i|>t\ia)\,\{1+o(1)\}\,,
$$
uniformly in $i$.  Theorem~3.3 follows from this property.

\bs


\baselineskip=16.7pt
\centerline{\bf REFERENCES}
\frenchspacing

\nh BENJAMINI, Y. AND HOCHBERG, Y. (1995).  Controlling the false
discovery rate:  a practical and powerful approach to multiple
testing.  {\sl J. Roy. Statist. Soc.} Ser.~B {\bf 57}, 289--300.

\nh BENJAMINI, Y. AND YEKUTIELI, D. (2001).   The control of the
    false discovery rate in multiple testing under dependency.  {\sl
    Ann. Statist.} {\bf 29}, 1165--1188.

\nh BENTKUS, V. AND G\"OTZE, F. (1996).  The Berry-Esseen bound for
Student's statistic.  {\sl Ann. Statist.} {\bf 24}, 491--503.

\nh BICKEL, P.J. and LEVINA, E. (2004).  Some theory for Fisher's
    linear discriminant, `naive Bayes', and some alternatives when there
    are many more variables than observations.  {\sl Bernoulli}, {\bf
    10}, 989-1010.

\nh DUDOIT, S., SHAFFER, J.P. AND BOLDRICK, J.C. (2003).
    Multiple hypothesis testing in microarray experiments.
    {\sl Statist. Sci.} {\bf 18}, 71--103.

\nh EFRON, B. (2004). Large-scale simultaneous hypothesis
        testing: the choice of a null hypothesis.
        {\sl J. Amer. Statist. Assoc.}
        {\bf 99}, 96--104.

\nh FAN, J., CHEN, Y., CHAN, H.M., TAM, P., AND REN, Y. (2005).
    Removing intensity effects and identifying significant
    genes for Affymetrix arrays in MIF-suppressed
    neuroblastoma cells.  {\sl Proc. Nat. Acad. Sci.  USA},
    {\bf 102}, 17751-17756.

\nh FAN, J. AND LI, R. (2006). Statistical Challenges with High
    Dimensionality: Feature Selection in Knowledge Discovery.
    Proceedings of the International Congress of Mathematicians (M.
    Sanz-Sole, J. Soria, J.L. Varona, J. Verdera, eds.) , Vol. III,
    595-622.

\nh  FAN, J., PENG, H., AND HUANG, T. (2005).  Semilinear
    high-dimensional model for normalization of microarray
    data:  a theoretical analysis and partial consistency  (With discussion).
    {\sl J. Amer. Statist. Assoc.} {\bf 100},
    781--813.

\nh  FAN, J., TAM, P., VANDE WOUDE, G. AND REN, Y. (2004).
    Normalization and analysis of cDNA micro-arrays
    using within-array replications applied to neuroblastoma
    cell response to a cytokine.  {\sl Proc. Nat. Acad. Sci. USA},
    {\bf 101}, 1135--1140.

\nh FINNER, H. AND ROTERS, M. (1998).  Asymptotic comparison of step-down and step-up multiple test procedures based on exchangeable test statistics. {\sl Ann. Statist.} {\bf 26}, 505--524.

\nh FINNER, H. AND ROTERS, M. (1999).  Asymptotic comparison of the critical values of step-down and step-up multiple comparison procedures. {\sl J. Statist. Plann. Inference} {\bf 79}, 11--30.

\nh FINNER, H. AND ROTERS, M. (2000).  On the critical value behaviour of multiple decision procedures. {\sl Scand. J. Statist.} {\bf 27}, 563--573.

\nh FINNER, H. AND ROTERS, M. (2001).  On the false discovery rate and expected type I errors. {\sl Biom. J.} {\bf 43}, 985--1005.

\nh FINNER, H. AND ROTERS, M. (2002).  Multiple hypotheses testing and expected number of type I errors. {\sl Ann. Statist.} {\bf 30}, 220--238.

\nh  GENOVESE, C. AND WASSERMAN, L. (2004).  A stochastic process
        approach to false discovery control. {\sl Ann. Statist.} {\bf
        32}, 1035--1061.


\nh HU, J. AND HE, X. (2006).   Enhanced quantile normalization of
        microarray data to reduce loss of information in the gene
        expression profile.  {\sl Biometrics}, to appear.

\nh HUANG, J., WANG, D., AND ZHANG, C. (2005). A two-way
    semi-linear model for normalization and significant analysis of
    cDNA microarray data. {\sl J. Amer. Statist. Assoc.} {\bf
    100}, 814--829.


\nh KOROSOK, M.R. AND MA, S. (2005).  Marginal asymptotics for the
    ``large $p$, small $n$'' paradigm:  With applications to micorarray
    data. {\sl Manuscript}.

\nh LEHMANN, E.L., ROMANO, J.P. AND SHAFFER, J.P. (2005).  On
        optimality of stepdown and stepup multiple test procedures.
        {\sl Ann. Statist.} {\bf 33}, 1084--1108.

\nh LEHMANN, E.L. AND ROMANO, J.P. (2005).  Generalizations of the
        familywise error rate. {\sl Ann. Statist.} {\bf 33}, 1138--1154.

\nh PETROV, V.V. (1975).  {\sl Sums of Independent Random
Variables.}  Springer, Berlin.


\nh  REINER, A., YEKUTIELI, D. AND BENJAMINI, Y. (2003).
        Identifying differentially expressed genes using false
        discovery rate controlling procedures. {\sl
        Bioinformatics}, {\bf 19}, 368-375.

\nh STOREY, J.D., TAYLOR, J.E., AND SIEGMUND, D. (2004).  Strong
control, conservative point estimation, and simultaneous
conservative consistency of
      false discovery rates: A unified approach. {\sl J. Roy. Stat. Soc.} Ser.~B
      {\bf 66}, 187--205


\nh VAN DER LAAN, M.J. AND BRYAN, J. (2001). Gene expression
    analysis with the parametric bootstrap. {\sl Biostatistics} {\bf
    2}, 445--461.

\nh WANG, Q. (2005).  Limit theorems for self-normalized large deviations.  {\sl Elect. J. Probab.} {\bf 10}, 1260--1285.

\vend